\documentstyle[12pt, amstex, amsfonts]{amsart}  

\title{Proof of the Hyperplane Zeros Conjecture of Lagarias and Wang}

\author{Wayne Lawton}

\address{Department of Mathematics\\
    National University of Singapore \\
    2 Science Drive 2\\
    Singapore 117543
    email: matwml(AT)nus.edu.sg}

\keywords {asymptotic, \'{e}tale mapping, expansive endomorphism
of torus group, Lojasiewicz's structure theorem for real analytic
sets, Pontryagin duality, resolution of singularities}

\subjclass {32C05, 37B05, 43A40}

\newtheorem{thm}{Theorem}[section]
\newtheorem{lem}{Lemma}[section]
\newtheorem{pro}{Proposition}[section]
\newtheorem{cor}{Corollary}[section]
\newtheorem{rem}{Remark}[section]
\newtheorem{defin}{Definition}[section]

\def\L2{L^2([ 0,2\pi))}
\def\L2R{L^2(\bold{R})}

\def\cA{{\cal A }}

\def\cC{{\cal C }}

\def\cE{{\cal E }}
\def\cF{{\cal F }}
\def\cG{{\cal G }}
\def\cH{{\cal H }}

\def\cM{{\cal M }}

\def\cS{{\cal S }}

\def\cU{{\cal U }}
\def\cV{{\cal V }}

\def\cZ{{\cal Z }}

\def\AA{{\Bbb A}}
\def\CC{{\Bbb C}}
\def\NN{{\Bbb N}}
\def\QQ{{\Bbb Q}}
\def\RR{{\Bbb R}}
\def\TT{{\Bbb T}}
\def\UU{{\Bbb U}}
\def\ZZ{{\Bbb Z}}

\begin{document} 

\begin{abstract}
We prove that a real analytic subset of a torus group that is
contained in its image under an expanding endomorphism is a finite
union of translates of closed subgroups. This confirms the
hyperplane zeros conjecture of Lagarias and Wang for real analytic
varieties. Our proof uses real analytic geometry, topological
dynamics and Fourier analysis.
\end{abstract}

\maketitle

\section{Introduction}
\setcounter{equation}{0}
By $\CC,$ $\RR,$ $\QQ,$ $\ZZ,$ $\NN \equiv \{1,2,3,4,...\}$ we
will denote the fields of complex, real, and rational numbers, the
ring of integers, and the set of natural numbers. For $n \in \NN,$
$\TT^n \equiv \RR^n/\ZZ^n$ denotes the $n$ dimensional torus
group, $\pi_n \, : \, \RR^n \rightarrow \TT^n$ denotes the
canonical homomorphism, and $\cE_n$ denotes the set of integer $n$
by $n$ expanding matrices (all eigenvalues have modulus $ > 1$).
If $E \in \cE_n,$ then $m \equiv |\det(E)| \geq 2,$ $E(\ZZ^n)$ is
a subgroup of $\ZZ^n$ with index $m,$ and $E$ induces an $m$ to
$1$ {\it expanding endomorphism} $E \, : \, \TT^n \rightarrow
\TT^n.$ For open $U \subset \RR^n,$ a function $h \, : \, U
\rightarrow \RR$ is called {\it real analytic} if $h \in
C^{\infty}(U)$ and the Taylor series of $h$ converges locally at
every point in $U.$ A function $h$ on $\RR^n$ is {\it periodic} if
$h(x + p) = h(x),  \, x \in \RR^n, \, p \in \ZZ^n.$ We denote the
zero set of a function $h$ by $\cZ_h.$ A subspace of $\RR^n$ is a
{\it rational subspace} if it is spanned by a finite subset of
$\ZZ^n.$ In \cite{Lagarias} Lagarias and Wang asserted the
following: \\ \\
{\it Hyperplane Zeros Conjecture : }
If $E \in \cE_n,$ $h : \RR^n \rightarrow \RR$ is real analytic and
periodic, and
$
    \cZ_h \subseteq E(\cZ_h) + \ZZ^n,
$
then
$
    \cZ_h = \bigcup_{i = 1}^{p} \, \left( V_i + x_i \right) + \ZZ^n,
$
where $p \geq 0$ and for every $i = 1,...,p,$ $V_i$ is a rational
subspace of $\RR^n$ and $x_i \in \RR^n.$
\begin{rem}
\label{rem:cerveau}
In \cite{cerveau} Cerveau, Conze and Raugi proved a result that is
similar to but weaker than the result asserted by the Hyperplane
Zeros Conjecture. In \cite{Lagarias} Lagarias and Wang used their
result to prove a result about tilings of $\RR^n$ that was
conjectured, and proved for the case $n = 1,$ by Gr\"{o}chenig and
Hass in \cite{grochenig}.
\end{rem}
Assume that $U \subseteq \TT^n$ is open. A function $f : U
\rightarrow \RR$ is called real analytic if $f \circ \pi_n$ is a
real analytic function on $\pi_{n}^{-1}(U),$ and $\cA(U)$ denotes
the ring of real analytic functions on $U.$ A subset $S \subseteq
U$ is called a {\it real analytic variety of} $U$ if there exists
$f \in \cA(U)$ such that $S = \cZ_f$ and is called a {\it real
analytic subset of} $U$ if for every $x \in U$ there exists an
open set $O_x \ni x$ and $f_x \in \cA(O_x)$ such that
$
    S \cap O_x = \{ \, y \, : \, y \in O_x, \, f_x(y) = 0 \, \}.
$
We let $\cV(U), \, \cS(U)$ denote the set of real analytic
varieties, subsets of $U.$ Clearly $\cV(U) \subseteq \cS(U).$
$\cG(\TT^n), \cG_c(\TT^n)$ denotes the set of closed, closed
connected subgroups of $\TT^n,$ and $\cF(\TT^n)$ denotes the set
of finite unions of translates of elements in $\cG_c(\TT^n).$ Then
$\cG_c(\TT^n) \subset \cG(\TT^n) \subset \cF(\TT^n)$ and Lemma
\ref{lem:FinV} implies that $\cF(\TT^n) \subset \cV(\TT^n).$
\begin{rem}
\label{rem:cartan} Cartan \cite{cartan} has constructed a compact
real analytic subset $S \subset \RR^3$ such that $f \in
\cA(\RR^3)$ and $f|_S = 0$ implies $f = 0.$ This implies that
$\cV(\TT^n) \neq \cS(\TT^n).$
\end{rem}
\begin{lem}
\label{lem:rational}
If $V$ is a subspace of $\RR^n,$ then $V$ is rational iff
$\pi_n(V) \in \cG_c(\TT^n).$
\end{lem}
\begin{pf}
This result was proved by Lagarias and Wang in \cite{Lagarias},
Theorem 4.1.
\end{pf}
Therefore the Hyperplane Zeros Conjecture is equivalent to the
statement of:
\begin{thm}
\label{thm:equiv}
If $E \in \cE_n,$ $S \in \cV(\TT^n),$ and $S \subseteq E(S),$ then
$S \in \cF(\TT^n).$
\end{thm}
The objective of this paper is to prove the following stronger
result:
\begin{thm}
\label{thm:main}
If $E \in \cE_n,$ $S \in \cS(\TT^n),$ and $S \subseteq E(S),$ then
$S \in \cF(\TT^n).$
\end{thm}
For $S \subseteq \TT^n,$ we let $\overline S$ denote the closure
of $S,$ $\partial S \equiv \overline S \, \backslash \, S$ denote
the boundary of $S,$ and $dim(S)$ denote the inductive dimension
of $S$ (Urysohn's Theorem implies that the small and large
inductive dimensions are equal since $\TT^n$ is separable and has
a countable basis). Then $-1 \leq dim(S) \leq n,$ and $dim(S) =
-1$ iff $S = \phi.$ For $x \in \TT^n$ we define
$
dim(S,x) \equiv \min \{ dim(S \cap O) \, : \, O \, \hbox{open},
\, x \in O \}.
$
Lojasiewicz'z Structure Theorem, described in Appendix C, shows
that if $S \in \cS(\TT^n)$ and $O \subseteq \TT^n$ is open, then
$O \cap S$ equals a union of submanifolds of $\TT^n,$ and
therefore $dim(O \cap S)$ equals the Lebesque covering dimension
of $O \cap S.$ Furthermore, if $S \in \cS(\TT^n),$ then $dim(S) =
n$ iff $S = \TT^n,$ and $dim(S) \leq 0$ iff $S$ is finite. These
facts justify our use of an induction procedure on the integers
$n$ and $dim(S)$ to prove Theorem \ref{thm:main}. For $B \subseteq
\TT^n$ we define
$
    \AA(B) \equiv \{ \, A \, : \, A \in \cS(\TT^n), \, B \subseteq A
    \, \}
$
and
$
    B^{\, *} \equiv \bigcap_{ \, A \in \AA(B) } A.
$
\begin{lem}
\label{lem:closure}
Assume that $E \in \cE_n.$ If $S \in \cS(\TT^n)$ then $E(S) \in
\cS(\TT^n).$ If $B \subseteq \TT^n$ and $B \subseteq E(B)$ then
$B^{\,*} \subseteq E(B^{\,*}).$
\end{lem}
\begin{pf}
The first assertion follows since $E$ is $m$ to $1$ and is locally
an analytic diffeomorphism. The second assertion then follows from
the definition of $B^{\, *}.$
\end{pf}
We call a subset $M$ of $\RR^n$ or of $\TT^n$ an $m$ {\it
dimensional real analytic submanifold} if it satisfies any of five
equivalent conditions described by Krantz in \cite{krantz}, p.
38-39. We let $\cM_m(\RR^n), \cM_m(\TT^n)$ denote the set of all
$m$ dimensional submanifolds of $\RR^n, \TT^n,$ and we let
$\cM(\RR^n), \cM(\TT^n)$ denote the set of all submanifolds of
$\RR^n, \TT^n,$ respectively. For $S \in \cS(\TT^n),$ a point $x
\in S$ is called a {\it regular point of dimension} $m$ if there
exists an open $O \subseteq \TT^n$ such that $S \cap O \in
\cM_m(\TT^n),$ and $R_m(S)$ denotes the set of regular points of
dimension.
\begin{lem}
\label{lem:density}
If $S \in \cS(\TT^n)$ and $d \equiv dim(S),$ then
$
    \overline {R_d(S)} = \{ x \in S \, : \, dim(S,x) = d \, \}.
$
\end{lem}
\begin{pf}
This follows from Theorem 1 in \cite{narasimhan}, p. 41.
\end{pf}
\begin{lem}
\label{lem:dimension}
If $S \in \cS(\TT^n)$ and $d \equiv dim(S),$ then
$\hbox{dim}(\, (S \, \backslash \, \, \overline {R_d(S)} \, \,
)^{\, *} \,) < d.$
\end{lem}
\begin{pf}
Since $A \equiv S \, \backslash \, R_d(S)$ is the set of singular
points in $S,$ Proposition 16 in \cite{whitney} implies that $A \in
\cS(\TT^n)$ and $\hbox{dim}(A) < d.$ Since $S \, \backslash
\overline \, {R_d(S)} \subseteq A$ it follows that
$A \in \AA\left(S \, \backslash \, \overline {R_d(S)}\right)$
hence $(S \, \backslash \, \, \overline {R_d(S)} \, \, )^{\, *}
\subseteq A.$ Therefore $\hbox{dim}(\, (S \, \backslash \, \,
\overline {R_d(S)} \, \, )^{\, *} \,) < d.$
\end{pf}
\begin{rem}
The result of Whitney-Bruhat in \cite{whitney} concerns real
analytic varieties in manifolds. Their result can be modified, by
considering intersections of $S$ with sufficiently small open
subsets $O \subseteq \TT^n$ such that $S \cap O \in \cV(O),$ to
obtain the result concerning real analytic subsets of $\TT^n$ that
we use in Lemma \ref{lem:dimension}.
\end{rem}
\begin{lem}
\label{lem:ERd}
If $E \in \cE_n,$ $S \in \cS(\TT^n),$ $d \equiv \hbox{dim}(S),$
and $E(S) \subseteq S,$ then it follows that $E(\overline
{R_d(S)}) \subseteq \overline {R_d(S)}.$
\end{lem}
\begin{pf}
Since $E$ is locally a diffeomorphism it preserves dimension.
\end{pf}
We now prove that Theorem \ref{thm:main} is equivalent to the
following:
\begin{thm}
\label{thm:redmain}
If $E \in \cE_n,$ $S\in \cS(\TT^n),$ $d \equiv \hbox{dim}(S),$ and
$S \subseteq E(S),$ then $S$ satisfies the following three
properties:
\begin{enumerate}
\item $S = E(S),$
\item $(S \, \backslash \, \, \overline {R_d(S)} \, \, )^{\, *} \in \cS(\TT^n),$
\item $\overline {R_d(S)} \in \cF(\TT^n).$
\end{enumerate}
\end{thm}
Assume that $E$ and $S$ satisfy the common hypothesis of Theorems
\ref{thm:main} and \ref{thm:redmain}. If $S \in \cF(\TT^n)$ then
$S$ satisfies the three properties in Theorem \ref{thm:redmain}.
Conversely, if $E$ and $S$ satisfy the three properties in Theorem
\ref{thm:redmain} then property (1) and Lemma \ref{lem:ERd} imply
that
$
    E(\overline {R_d(S)}) \subseteq \overline {R_d(S)}.
$
Therefore
$S \, \backslash \, \, \overline {R_d(S)} \subset E(S \,
\backslash \, \, \overline {R_d(S)})$
and Lemma \ref{lem:closure} implies that
$(S \, \backslash \, \overline {R_d(S)})^* \subseteq E((S \,
\backslash \, \overline {R_d(S)})^*).$
This inclusion, together with property (2), implies that $(S \,
\backslash \, \overline {R_d(S)})^*$ satisfies the hypothesis of
Theorem \ref{thm:main}, and Lemma \ref{lem:dimension} implies that
it has dimension $< dim(S),$ hence by induction $(S \, \backslash
\, \overline {R_d(S)})^* \in \cF(\TT^n).$ This inclusion and
property (3) implies that $S \in \cF(\TT^n)$ and completes the
proof.
\\ \\
The remaining four sections of this paper show that if $S$ satisfies
the hypothesis of Theorem \ref{thm:redmain}, then $S$ satisfies
properties (1), (2), and (3). Section 2 uses stationary properties
of real analytic sets to show properties (1) and (2). Sections 3, 4,
and 5 use topological dynamics of mappings associated with $E,$
together with the induction hypothesis on $n$ and $\hbox{dim}(S),$
to show property (3). Section 3 derives an asymptotic property of
the map $E \, : \, \RR^n \, \rightarrow \RR^n$ on submanifolds of
$\RR^n$ and uses it to derive a sufficient condition to ensure that
$S$ satisfies the following invariance property: there exists $y \in
S$ and $H \in \cG_c(\TT^n)$ such that $\hbox{dim}(H) > 1$ and $y + H
\subseteq S.$ Section 4 uses results about the Hausdorff topology on
$\cG_c(\TT^n),$ derived in Appendix B, to prove that $S$ satisfies
property (3) whenever it satisfies special invariance properties.
Section 5 uses Lojaciewicz'z Structure Theorem for Real Analytic
Sets, stated in Appendix C, together with recent results of Hiraide
about the topological dynamics of positively expansive maps, to
construct a resolution of singularities for the set $\overline
{R_d(S)}.$ It uses this resolution to prove that $S$ satisfies these
special invariance properties that ensure that it satisfies property
(3). This completes the proof of Theorems \ref{thm:redmain} and
\ref{thm:main} and confirms the Hyperplane Zeros Conjecture.
\section{Stationarity Properties of Real Analytic Sets}
\setcounter{equation}{0}
Let $X$ be a set and $K \subseteq X.$ A family $S_\alpha \subseteq
X, \, \alpha \in I$ indexed by a partially ordered set $(I,\leq)$ is
called a {\it decreasing filtered family} (DFF) of subsets of $X$ if
$S_\beta \subseteq S_\alpha, \, \alpha \leq \beta.$ It is {\it
stationary on} $K$ if there exists $\alpha \in I$ such that $S_\beta
\cap K = S_\alpha \cap K, \, \alpha \leq \beta.$
\begin{pro}
\label{pro:stat}
The following assertions hold:
\begin{enumerate}
\item If $S_\alpha \in \cS(\RR^n), \, \alpha \in I$ is a
DFF, then it is stationary on every compact subset.
\item If $S_\alpha \in \cS(\TT^n), \, \alpha \in I$ is a
DFF, then it is stationary on $\TT^n.$
\item If $S_\alpha \in \cS(\TT^n), \alpha \in I$ is
an arbitrary family of subsets, then $\bigcap_{\alpha \in I}
S_\alpha \in \cS(\TT^n).$
\end{enumerate}
\end{pro}
\begin{pf}
Narasimhan proved the first assertion in \cite{narasimhan},
Corollary 1, p. 99. To show the second assertion construct a compact
$K \subset \RR^n$ such that $\pi_n(K) = \TT^n.$ Then
$\pi_{n}^{-1}(S_\alpha) \in \cS(\RR^n)$ is a DFF, hence it is
stationary on $K,$ hence $S_\alpha$ is stationary on $\TT^n.$ An
alternative proof can be based on the fact that $\TT^n$ is a compact
real analytic subset of $\RR^{2n}.$ The third assertion follows by
applying the second assertion to the DFF consisting of finite
intersections of elements in $\{\, S_\alpha \, : \, \alpha \in I \,
\}$ (indexed by the partially ordered set of finite subsets of $I$).
\end{pf}
\begin{cor}
\label{cor:equal}
If $E \in \cE_n,$ $S \in \cS(\TT^n),$ and $S \subseteq E(S),$ then
$S = E(S).$
\end{cor}
\begin{pf}
If $S \neq E(S)$ then the sequence defined by $S_p := \bigcap_{\, k
= 0}^{\, p} \, E^{-k}(S), \ p \in \NN$ is a strictly decreasing
filtered family of subsets of $\cS(\TT^n)$ (indexed by the set $\NN$
with standard partial order) thus contradicting Proposition
\ref{pro:stat}.
\end{pf}
\begin{cor}
\label{cor:closure}
For every $B \subseteq \TT^n,$ $B^* \in \cS(\TT^n).$
\end{cor}
\begin{pf}
Follows by applying assertion (3) of Proposition \ref{pro:stat} to
the family $\cC(B).$
\end{pf}
\begin{rem}
\label{rem:frisch}
Frisch \cite{frisch}, Theorem I,9 proved that the ring $\cA(\TT^n)$
is Noetherian. Therefore, by a standard result \cite{hart}, if
$S_\alpha \in \cV(\TT^n), \, \alpha \in I$ is a decreasing filtered
family, then it is stationary on $\TT^n.$ Furthermore, if $S_\alpha
\in \cV(\TT^n), \alpha \in I$ is an arbitrary family, then
$\bigcap_{\alpha \in I} S_\alpha \in \cV(\TT^n).$ Also see
\cite{kiehl}.
\end{rem}
\section{Asymptotic Tangent Vectors}
\setcounter{equation}{0}
Let $(X,\rho)$ be a metric space. For $A, B \subseteq X$ we define
the asymmetric distance from $A$ to $B$ by
$$
    \rho(A,B) := \sup_{ \, a \, \in \, A \, } \
    \inf_{ \, b \, \in \, B \, } \ \rho(a,b).
$$
For $r > 0$ and $x \in X$ we define the closed unit ball of radius
$r> 0$ centered at $x$ by
$$
    B_{\rho}(r,x) := \{ \, y \in X \, : \, \rho(x,y) \leq r \, \}.
$$
For $x, v \in \RR^n$ with $v \neq 0$ we define the line through
$x$ in the direction of $v$ by
$$
    \ell(x,v) := x + \hbox{span}(\{ v \}).
$$
If $k \geq 1,$ $M$ is a $k$ dimensional $C^{1}$ submanifold of
$\RR^n,$ and $x \in M,$ then we identity the tangent space
$T_x(M)$ to $M$ at $x$ with a $k$ dimensional subspace of $\RR^n$
and we observe that the set $x + T_x(M)$ is a $k$ dimensional
affine subset of $\RR^n$ that is geometrically tangent to $M$ are
$x.$
\begin{lem}
\label{lem:taylor1}
Let $|| \, || \, : \RR^n \rightarrow (0,\infty)$ be a norm on
$\RR^n$ and let $\rho$ be the associated metric on $\RR^n$ defined
by $\rho(x,y) = ||x-y||, \, x, y \in \RR^n.$ If $M$ is a $C^{2}$
submanifold of $\RR^n$ then there exists a continuous function
$\alpha \, : \, M \rightarrow (0,\infty)$ such that
\begin{equation}
\label{eqn:taylor2}
   \rho(B_{\rho}(r,x) \cap (x + T_x(M)),M) \leq \alpha(x) \, r^2,
   \ \ x \in M, \ r \in [0,1].
\end{equation}
\end{lem}
\begin{pf}
Follows from the error bound for the first degree Taylor
approximation.
\end{pf}
Throughout the remainder of this section we assume that $E \in
\cE_n.$ If $j \geq 0$ then $E^j(\ell(x,v)) = \ell(E^j x, E^j v),$
$E^j(M)$ is a $k$ dimensional $C^{1}$ submanifold of $\RR^n,$ and
$T_{E^jx}(E^j(M)) = E^j(T_x(M)).$ If $|| \, ||$ is a norm on
$\RR^n$ and $F$ is an $n$ by $n$ matrix then we define the
associated matrix norm by
$
    ||F|| \equiv \max \{ \, ||F v|| \, : \, v \in \RR^n, \, ||v||
    = 1 \, \}.
$
\begin{defin}
\label{def:asymptotic}
A triplet $(M,x,v),$ where $M \in \cM(\RR^n),$ $x \in M,$ $v \in
T_x(M),$ and $v \neq 0,$ is $E$ asymptotic if there exists a norm
$|| \, ||$ on $\RR^n$ and an associated metric $\rho$ on $\RR^n$
such that
$
    \lim_{j \rightarrow \infty}
    \rho(\, B_{\rho}(1, E^j x) \cap \ell(E^j x,E^j v)), \, E^j(M)) = 0.
$
\end{defin}
\noindent We observe that this concept is independent of the norm.
\begin{lem}
\label{lem:flat}
If $M \in \cM(\RR^n)$ then there exists a continuous function
\newline $\alpha \, : \, M \rightarrow (0,\infty)$ such that for $x \in
M,$ $v \in T_x(M),$ $v \neq 0,$
\begin{equation}
\label{eqn:taylor3}
\rho(\, B_{\rho}(1, E^j x) \cap \ell(E^j x,E^j v)), \, E^j(M))
\leq \alpha(x) \, ||v||^2 \, \frac{||E^j||}{||E^j v||^{2}}.
\end{equation}
\end{lem}
\begin{pf}
Clearly
$B_{\rho}(1, E^j x) \cap \ell(E^j x,E^j v)) \subseteq
E^j(B_{\rho}(r_j, x) \cap \ell(x,v)), \, r_j \equiv ||v|| \, ||E^j
v||^{-1}.$
Therefore
$\rho(\, B_{\rho}(1, E^j x) \cap \ell(E^j x,E^j v)), E^j(M)) \leq
\rho(\, E^j(B_{\rho}(r_j, x) \cap \ell(x,v)), \, E^j(M))$
hence
$\rho(E^j(B_{\rho}(r_j, x) \cap \ell(x,v)), \, E^j(M)) \leq
||E^j|| \, \rho(B_{\rho}(r_j, x) \cap \ell(x,v), \, M).$
We observe that since $\ell(x,v) \subseteq x + T_x(M),$ Lemma
\ref{lem:taylor1} implies that there exists a continuous function
$\alpha \, : \, M \rightarrow (0,\infty)$ such that
$\rho(B_{\rho}(r_j, x) \cap \ell(x,v),\,M) \leq \alpha(x)\,
r_{j}^{2}$
and Inequality \ref{eqn:taylor3} follows from combining these
three inequalities above.
\end{pf}
Let $\Lambda(E)$ denote the set of eigenvalues of $E,$ let
$\sigma := \max \{ \, | \lambda | \, : \, \lambda \in
\Lambda(E)\}$
denote the spectral radius of $E,$ and let
$\sigma_1 := \max \{ \, | \mu | \, : \, \mu \in \Lambda(E) \,
\hbox{and} \, |\mu| < \sigma \, \}.$
Let
$
    V_{\lambda} := \{ \, x \in V \, : \, (E-\lambda)^n x = 0 \}
$
denote the $E$ invariant subspace associated to $\lambda \in
\Lambda(E).$ Define subspaces
$
    V_{\sigma} := \sum_{|\lambda| = \sigma} \, V_{\lambda}
$
and
$
    V_{\sigma}^{\perp} := \sum_{|\mu| < \sigma} \, V_{\mu}.
$
\begin{thm}
\label{thm:asymptotic}
If $M \in \cM(\RR^n),$ $x \in M,$ $v \in T_x(M),$ and $v \notin
V_{\sigma}^{\perp},$ then the triplet $(M,x,v)$ is $E$ asymptotic.
\end{thm}
\begin{pf}
$\lim_{j \rightarrow \infty} \frac{||E^j||}{||E^j v||^{2}} = 0$
since $||E^j|| \approx \sigma^j$ and $||E^j v|| \geq \beta \sigma^j$
for some $\beta > 0.$
\end{pf}
The following result illustrates how asymptotic properties of the
$E$ imply invariance properties of $S \in \cS(\TT^n)$ that satisfy
$E(S) = S.$
\begin{pro}
\label{prop:asyminv}
If $E \in \cE_n,$ $S \in \cS(\TT^n),$ $E(S) = S,$ $M \in
\cM(\RR^n),$ $\pi_n(M) \subseteq S,$ $x \in M,$ then at least one of
the following properties hold:
\begin{enumerate}
\item $T_x(M) \subset V_{\sigma}^{\perp},$
\item there exists $y \in S$ and $H \in \cG_c(\TT^n)$
        such that $\hbox{dim}(H) > 1$ and $y + H \subseteq S.$
\end{enumerate}
\end{pro}
\begin{pf}
Assume that $v \in T_x(M)$ and $v \notin V_{\sigma}^{\perp}.$ It
suffices to show that $S$ satisfies property (2) above. We choose a
norm $|| \, ||$ on $\RR^n$ and let $\rho$ denote the corresonding
metric on $\RR^n.$ We also let $\rho$ denote the metric on $\TT^n$
defined by $\rho(a,b) \equiv \min \{ \rho(\alpha,\beta) \, : \,
\pi_n(\alpha) = a, \, \pi_n(\beta) = b \, \}.$ Since $\TT^n$ and the
set of unit vectors in $\TT^n$ are compact there exists a function
$j \, : \, \NN \rightarrow \NN$ such that satisfies the following
properties
\begin{enumerate}
\item $\lim_{i \rightarrow \infty} j(i) = \infty,$
\item there exists $y \in S$ such that the sequence $y_i \equiv \pi_n(x_i),$
where $x_i \equiv E^{j(i)}x,$ converges to $y,$
\item there exists $u \in \RR^n$ such that the sequence
$u_i \equiv \frac{E^{j(i)}v}{||E^{j(i)}v||}$ converges to $u.$
\end{enumerate}
We construct the vector spaces $U_i \equiv \hbox{span} \{ \, u_i \,
\}, \, i \in \NN$ and $U \equiv \hbox{span} \{ \, u \, \}$ and
construct $H \equiv \overline {\pi_n(U)}.$ Clearly $H \in
\cG_c(\TT^n)$ and $\hbox{dim}(H) \geq 1.$ It suffices to show that
$y + H \subseteq S.$ We construct the sequence of submanifolds $M_i
\equiv E^{j(i)}(M), \, i \in \NN$ of $\RR^n.$ Clearly $\pi_{n}(M_i)
\subseteq S, \, i \in \NN$ since $\pi_n(M) \subseteq S$ and $E(S) =
S.$ Theorem \ref{thm:asymptotic} implies that $(M,x,v)$ is an
asymptotic triple in the sense of Definition \ref{def:asymptotic}.
Therefore the sequence $\rho(\, x_i + B_{\rho}(1,0) \cap U_i, \, M_i
\,)$ converges to zero. Therefore the sequence
$\rho(\, \pi_n(x_i + B_{\rho}(1,0) \cap U_i), \, \pi_n(M_i) \,)$
converges to zero, hence the sequence
$\rho(\, \pi_n(x_i + B_{\rho}(1,0) \cap U_i), \, S \,)$
converges to zero. Since $\pi_n(x_i) = y_i$ converges to $y$ and
$B_{\rho}(1,0) \cap U_i$ converges to $B_{\rho}(1,0) \cap U,$ and
$S$ is compact,
$\rho(\, y + \pi_n(B_{\rho}(1,0) \cap U), \, S \,) = 0.$
Therefore $y + \pi_n(B_{\rho}(1,0) \cap U) \subseteq S.$ Since $S$
is a real analytic set $y + \pi_n(U) \subseteq S.$ Since $S$ is
closed $y + H \subseteq S.$
\end{pf}
\section{Invariance Properties of Subsets of $\TT^n$}
\setcounter{equation}{0}
We derive properties of certain subsets of $\TT^n$ that are either
invariant under translation by elements in $G,$ where $G \in
\cG_c(\TT^n),$ or that are subsets of $\pi_n(V),$ where $V$ is a
proper subspace of $\RR^n$ that is invariant under $E \in \cE_n.$
\begin{defin}
\label{defin:invariant}
For $S \subseteq \TT^n$ and $G \in \cG_c(\TT^n)$ we define $S_G
\equiv \{ x \in S \, : \, x+G \subseteq S \, \}.$ The set $S_G$ is
called the $G$ {\it invariant subset of} $S.$
\end{defin}
\begin{lem}
\label{lem:inv1}
If $S \in \cS(\TT^n)$ and $G \in \cG_c(\TT^n),$ then $S_G \in
\cS(\TT^n).$
\end{lem}
\begin{pf}
Follows from assertion (3) of Proposition \ref{pro:stat} since
$S_G = \bigcap_{\,g \in G} \left(S - g\right).$
\end{pf}
For $G \in \cG_c(\TT^n)$ we let $\pi_G \, : \, \TT^n \rightarrow
\TT^n \, / \, G$ denote the canonical homomorphism. We observe
that $\TT^n \, / \, G$ is isomorphic to $\TT^m$ where $m = n -
dim(G).$ Therefore the sets $\cS(\TT^n \, / \, G)$ and $\cF(\TT^n
\, / \, G)$ are defined. Corollary \ref{cor:GtoH} implies that
there exists $H \in \cG_c(\TT^n)$ such that $G \cap H = \{0\}$ and
$G + H = \TT^n.$ We further observe that $\pi_G|_{H} \, : \, H
\rightarrow \TT^n \, / \, G$ is an isomorphism and that $\pi_G|_H$
induces bijections between $\cS(H)$ and $\cS(\TT^n \, / \, G)$ and
between $\cF(H)$ and $\cF(\TT^n \, / \, G).$
\begin{lem}
\label{lem:inv2}
If $S \in \cS(\TT^n)$ and $G \in \cG_c(\TT^n),$ then $\pi_G(S_G)$
satisfies:
\begin{enumerate}
\item $dim(\pi_G(S_G)) = dim(S_G) - dim(G),$
\item $S_G \in \cS(\TT^n)$ iff $\pi_G(S_G) \in \cS(\TT^n \, / \, G),$
\item $S_G \in \cF(\TT^n)$ iff $\pi_G(S_G) \in \cF(\TT^n \, / \, G),$
\item if $E \in \cE_n$ and $E(G) = G,$ then $E$ induces an expanding endomorphism
\newline $E \, : \, \TT^n \, / \, G \rightarrow \TT^n \, / \, G,$ and if
$S \subseteq E(S)$ then $\pi_G(S_G) \subseteq E(\pi_G(S_G)).$
\end{enumerate}
\end{lem}
\begin{pf}
If $H \in \cG_c(\TT^n)$ is the subgroup in Corollary
\ref{cor:GtoH}, then $S_G = (H \cap S_G) + G$ and $H \cap G =
\{0\}$ hence $dim(S_G) = dim(H \cap S_G) + dim(G).$ Since
$\pi_G|_H \, : \, H \rightarrow \TT^n \, / \, G$ is an analytic
bijection, $dim(\pi_G(H \cap S_G)) = dim(H \cap S_G)$ . Property
(1) follows since $\pi_G(S_G) = \pi_{G}(H \cap S_G).$ If $S_G \in
\cS(\TT^n),$ then $H \cap S_G \in \cS(H).$ Therefore, since
$\pi_G|H \, : \, H \rightarrow \TT^n \, / \, G$ induces a
bijection between $\cS(H)$ and $\cS(\TT^n \, / \, G),$ it follows
that $\pi_G(S_G) = \pi_G(H \cap S_G) \in \cS(\TT^n \, / \, G).$
Since $\pi_G \, : \, \TT^n \rightarrow \TT^n \, / \, G$ is real
analytic and $S_G = \pi_{G}^{-1}(\pi_G(S_G),$ it follows that if
$\pi_G(S_G) \in \cS(\TT^n \, / \, G)$ then  $\pi_G(S_G) \in
\cS(\TT^n).$ This proves property (2). Properties (3) and (4) are
evident.
\end{pf}
\begin{lem}
\label{lem:inv3}
Assume that $E \in \cE_n,$ $S \in \cS(\TT^n),$ and $S \subseteq
E(S),$ and there exists $H \in \cG_c(\TT^n),$ $dim(H) \geq 1,$ and
$S = S_H.$ Then under the induction hypothesis on $dim(S),$ $S \in
\cF(\TT^n).$
\end{lem}
\begin{pf}
Theorem \ref{thm:invinv} implies that there exists $G \in
\cG_c(\TT^n)$ and $p \in \NN$ such that $dim(G) \geq 1,$ $E^p(G) =
G,$ and $S = S_G.$ Clearly $E^p \in \cE_n$ and Corollary
\ref{cor:equal} implies that $E^p(S) = S.$ Lemma \ref{lem:inv2}
implies that $\pi_G(S_G) \in \cS(\TT^n \, / \, G),$ $\pi_G(S_G)
\subseteq E^p(\pi_G(S_G)),$ and $dim(\pi_G(S_G)) = d - dim(G) < d.$
Therefore by induction on $d,$ $\pi_G(S_G) \in \cF(\TT^n \, / \, G)$
hence Lemma \ref{lem:inv2} implies that $S_G \in \cF(\TT^n).$
\end{pf}
\begin{lem}
\label{lem:pinxy}
If $E \in \cE_n,$ $x \in \TT^n,$ and $V \subseteq \RR^n$ is a
subspace that satisfies \newline $E(\pi_n(V) + x) = \pi_n(V) + x,$
then $E(\pi_n(V)) = \pi_n(V),$ $E(V) = V,$ and there exists $y \in
\TT^n$ such that $E(y) = y$ and $\pi_n(V) + y = \pi_n(V) + x.$
\end{lem}
\begin{pf}
Since $E(\pi_n(V) + x) = E(\pi_n(V)) + E(x)$ the assumption
$E(\pi_n(V) + x) = \pi_n(V) + x$ implies that $E(\pi_n(V)) =
\pi_n(V) + x - E(x).$ Therefore, since both $\pi_n(V)$ and
$E(\pi_n(V))$ are subgroups of $\TT^n,$ $E(\pi_n(V)) = \pi_n(V),$
$E(V) = V,$ and $x - E(x) \in \pi_n(V).$ Let $I \, : \, \TT^n
\rightarrow \TT^n$ denote the identity map. Since $E \in \cE_n$
the map \newline $(I - E)|_{\pi_n(V)} \, : \, \pi_n(V) \rightarrow
\pi_n(V)$ is surjective. Therefore there exists $z \in \pi_n(V)$
such that $z - E(z) = -x + E(x).$ Let $y \equiv x + z.$ Then $E(y)
= E(x + z) = x+z = y$ and $\pi_n(V) + y = \pi_n(V) + x.$
\end{pf}
For every subspace $V \subseteq \RR^n$ we let $V_{rat}$ denote the
subspace spanned by $V \cap \ZZ^n.$ Clearly $V_{rat}$ is the
largest rational subspace contained in $V$ and $\pi_n(V_{rat}) \in
\cG_c(\TT^n).$
\begin{lem}
\label{lem:EHE}
If $E \in \cE_n$ and $V \subseteq \RR^n$ is a subspace that
satisfies $E(V) = V,$ then $E(\pi_n(V_{rat})) = \pi_n(V_{rat}).$
\end{lem}
\begin{pf}
Since $E(V_{rat}) = \hbox{span }(V \cap E(\ZZ^n)) \subseteq
V_{rat},$ and $E \, : \, \RR^n \rightarrow \RR^n$ is injective,
and $V_{rat}$ is finite dimensional, $E(V_{rat}) = V_{rat}.$
Therefore $E(\pi_n(V_{rat})) = \pi_n(V_{rat}).$
\end{pf}
In the following result we assume that $\RR^n,$ $\TT^n,$ and certain related
quotient groups are equipped with a Riemannian structure defined by the standard Euclidean scalar product
on $\RR^n.$ The length of a $C^1$ parameterized path $\gamma \, : \, [0,1] \rightarrow Y,$
where $Y$ is any Riemannian manifold, is defined by $\int_0^1 ||\frac{d\gamma}{dt}||\,dt.$
We observe that if $\gamma$ is a $C^1$ path in $\TT^n,$ then $\gamma \circ \pi_n$ is
a path in $\RR^n$ and that both paths have the same lengths. We observe that if $m \in \NN,$
$W$ is a subspace of $\RR^m,$ $\psi \, : \, W \rightarrow \TT^m$ is a injective homomorphism
such that $||\, d \psi (v)|| = ||v||, \, v \in W,$ and $\gamma$ is a $C^1$ path in
$\psi(W),$ then there exists a unique $C^1$ path $\gamma_W$ in $W$ such that $\gamma = \psi \circ \gamma_W$
and that the length of $\gamma_W$ equals the length of $\gamma.$ It follows that if $K \subset \pi_n(W),$
$c > 0,$ and every two points in $K$ are connected by a $C^1$ path having length $\leq c,$
then $\psi^{-1}(K)$ is a bounded subset of $W.$
\begin{pro}
\label{pro:irr}
If $E \in \cE_n,$ $x \in \TT^n,$ $V \subseteq \RR^n$ is a subspace
that satisfies $E(V) = V,$ $K \subseteq \pi_n(V) + x$ is a
nonempty, closed subset that satisfies $E(K) = K,$ $c > 0,$ and every two points in
$K$ are connected by a $C^1$ path having length $\leq c,$ then there
exists $y \in \TT^n$ such that $E(y) = y$ and $K \subseteq
\pi_n(V_{rat}) + y.$
\end{pro}
\begin{pf}
Clearly $K \subseteq \pi_n(V) + E(x),$ therefore, since $\pi_n(V)$
is a subgroup of $\TT^n$ and since $K$ is nonempty, $\pi_n(V) +
E(x) = \pi_n(V) + x$ hence $E(\pi_n(V) + x) = \pi_n(V) + x.$
Therefore Lemma \ref{lem:pinxy} implies that there exists $y \in
\TT^n$ such that $E(y) = y$ and $\pi_n(V) + y = \pi_n(V) + x.$
Construct the set $J \equiv K - y.$ Then $J$ is closed, $J
\subseteq \pi_n(V),$ and $E(J) = J.$ It suffices to prove that $J
\subseteq \pi_n(V_{rat}).$ We construct the quotient groups
$W \equiv V \, / \, V_{rat}, $
and
$G \equiv  \pi_n(V) \, / \, \pi_n(V_{rat}),$
and we let $\psi \, : \, W \rightarrow G$ denote the unique
homomorphism that makes the following diagram commute:
$$
\begin{array}{ccccccc}
V_{rat} & \longrightarrow & V & \pi_W \atop \longrightarrow & W & \longrightarrow & \RR^n \, / \,  V_{rat} \\ \\
\downarrow  \pi_n &       & \downarrow \pi_n  &     & \downarrow \psi &  & \downarrow \pi_n \\ \\
\pi_{n}(V_{rat}) & \longrightarrow & \pi_n(V) & \pi_G \atop \longrightarrow & G & \longrightarrow & \TT^n \, / \, \pi_n(V_{rat})\\ \\
\end{array}
$$
Here all unlabelled right arrows denote inclusion maps, and $\pi_n,$
$\pi_W$ and $\pi_G$ denote canonical epimorphisms. Furthermore, we
observe that since the horizontal sequences are exact and since
the kernel of $\pi_{n} \, : \, V \rightarrow \pi_{n}(V)$ equals $V
\cap \ZZ^n \subset V_{rat},$ the map $\psi$ is injective. We
construct the set
${\widetilde J} \equiv \psi^{-1}(\pi_G(J)).$
We observe that since $K$ has the property that every two points in $K$
can be connected by a $C^1$ path having length $\leq c,$ then both $J$ and
$\pi_G(J)$ have this property. It follows from the preceding discussion
that ${\widetilde J}$ is bounded. Since $E(V_{rat}) = V_{rat},$ $E$
induces an expansive endomorphism $E \, : \, W \rightarrow W.$
Since the diagram commutes and $E(J) = J,$ it follows that
$E(\widetilde J) = \widetilde J.$ Therefore, since $\lim_{j\rightarrow \infty} ||E^jw|| = \infty$
for every $w \in W \backslash \{0\},$ and since ${\widetilde J}$ is bounded, it follows that
$\widetilde J = \{0\}.$ Therefore $\pi_G(J) = \{0\}$ hence $J \subseteq \pi_{n}(V_{rat}).$
\end{pf}
\section{\'{E}tale Construction and Resolution of Singularities}
\setcounter{equation}{0}
If $X$ is a topological space and $x \in X,$ we introduce an
equivalence relation $\approx_{x}$ on the set of subsets of $X$ as
follows: for subsets $M, N \subseteq X,$
$
    M \approx_x N
$
if there exists an open neighborhood $O$ of $x$ such that
$
    M \cap O = N \cap O \, ,
$
and we denote the corresponding equivalence class of a set $M$ by
$M_x.$ The set $M_x$ is called the germ of the set $M$ at $x.$
Clearly $\phi \in M_x$ if and only if $x \notin M,$ and $\{x\} \in
M_x$ if and only if $x$ is an isolated point in $M.$ A mapping of
topological spaces $f : X \rightarrow Y$ is called an \'{e}tale
mapping if every point $x \in X$ has a neighborhood $O$ such that
the restriction $f|_{O} \, : \, O \rightarrow f(O)$ is a
homeomorphism. Clearly \'{e}tale mappings are continuous and open.
We refer the reader to Godement \cite{godement} for a more
detailed discussion of germs and \'{e}tale maps. We fix $S \in
\cS(\TT^n)$ and let $d = dim(S).$ For $x \in \overline {R_d(S)},$
we let $[x]$ denote the set of germs $M_x$ where $M \in
\cM_d(\TT^n),$ $M \subseteq \overline {R_d(S)},$ and $x \in M.$ We
define $\widetilde S_0 \equiv \bigcup_{x \in \overline {R_d(S)}} \
[x],$ and construct the map $\tau \, : \, \widetilde S_0
\rightarrow \overline {R_d(S)}$ so that $\tau(y) = x$ if and only
if $y \in [x].$ If $x \in R_d(S),$ then $[x]$ has one point,
$R_d(S) \subseteq \tau(\widetilde S_0) \subseteq
\overline{R_d(S)},$ and the restriction $\tau|_{\tau^{-1}(R_d(S))}
\, : \, \tau^{-1}(R_d(S)) \rightarrow R_d(S)$ is a bijection.
We construct a topology on $\widetilde S_0$ generated by the
(open) sets $\{ \, M_y \, : \, y \in M \, \}$ where $M \in
\cM_d(\TT^n)$ and $M \subseteq \overline {R_d(S)}.$ We let
$\cC(\widetilde S_0)$ denote the set of connected components of
$\widetilde S_0.$
\begin{lem}
\label{lem:etal}
The space $\widetilde S_0$ is Hausdorff, $\tau \, : \, \widetilde
S_0 \rightarrow \tau(\widetilde S_0)$ is an \'{e}tale mapping,
finite to one, and proper (the inverse image of every compact set
is compact) and its restriction $\tau|_{\tau^{-1}(R_d(S))} \, : \,
\tau^{-1}(R_d(S)) \rightarrow R_d(S)$ is a homeomorphism. The set
$\cC(\widetilde S_0)$ is finite.
\end{lem}
\begin{pf}
These first assertions follow from Proposition \ref{pro:loja} and
Corollary \ref{cor:loja1}.
\end{pf}
Since $\tau$ maps $\widetilde S_0$ locally onto analytic
submanifolds of $\TT^n,$ it induces the structure of a real
analytic manifold on $\widetilde S_0$ and then $\tau \, : \, S_0
\rightarrow \TT^n$ is a real analytic immersion.
If $M \in \cM(\TT^n)$ and $x \in M$ we let $T_x(M)$ denote the
tangent space to $M$ at $x.$ We identify every tangent space to
$\TT^n$ with $\RR^n$ (via the canonical homomorphism $\pi_n \, :
\, \RR^n \rightarrow \TT^n$) and we identify $T_x(M)$ with an
$dim(M)$-dimensional subspace of $\RR^n.$ If $S \in \cS(\TT^n)$
and $x$ is a regular point in $S$ then $T_x(S)$ denotes the
tangent space to $S$ at $x.$
We recall that a Riemannian structure on a differentiable manifold
consists of a symmetric positive definite bilinear form at each of
its tangent spaces, and that a Riemannian structure defines a
metric space whose distance function is the corresponding geodesic
distance. We consider $\TT^n$ to have the Riemannian structure
given by the standard bilinear form on $\RR^n.$
Since $\tau$ is a smooth \'{e}tale mapping, it induces a
Riemannian structure on $\widetilde S_0$ such that $\tau$ maps
each tangent space of $\widetilde S_0$ isometrically into $\RR^n.$
We consider $\widetilde S_0$ to be a metric space (defined by this
Riemannian structure) and we let $\widetilde S$ denote the
completion of the metric space $\widetilde S_0.$ Since $\tau$ is
proper, $\widetilde S$ is compact. Since $\tau \, : \, \widetilde
S_0 \rightarrow \TT^n$ is uniformly continuous, $\tau$ extends to
define a continuous surjection $\tau \, : \, \widetilde S
\rightarrow \overline {R_d(S)}.$ An analysis based on Proposition
\ref{pro:loja} shows that $\widetilde S$ is locally connected and
has a finite number of connected components.
\begin{lem}
\label{lem:lift}
If $E \in \cE_n,$ $S \in \cS(\TT^n),$ and $E(S) = S,$ then
$E(\tau(\widetilde S_0)) \subseteq \tau(\widetilde S_0)$ and there
exists a unique map $\widetilde E \, : \, \widetilde S_0
\rightarrow \widetilde S_0$ that satisfies $\tau \circ \widetilde
E = E \circ \tau \, : \, \widetilde S_0 \rightarrow \widetilde
S_0.$ The map $\widetilde E \, : \, \widetilde S_0 \rightarrow
\widetilde S_0$ is surjective, locally an analytic diffeomorphism,
and it extends to a map $\widetilde E \, : \widetilde S
\rightarrow \widetilde S$ that satisfies $\tau \circ \widetilde E
= \tau \circ E \, : \, \widetilde S \rightarrow \overline
{R_d(S)}.$ $\widetilde E$ maps each component of $\widetilde S_0$
onto a component of $\widetilde S_0$ by a locally analytic
diffeomorphism and induces a permutation of $\cC(\widetilde S_0)$
hence there exists $p \in \NN$ such that
$\widetilde E^p(C) = C, \ C \in \cC(\widetilde S_0).$
For each $C \in \cC(\widetilde S_0),$ the closure $\overline C$ is
compact, connected, and locally connected and $\widetilde E^p$
extends to give a positively expansive map continuous surjection
$\widetilde E^p \, : \, \overline C \rightarrow \overline C.$
\end{lem}
\begin{pf}
Assume that $x \in \tau(\widetilde S_0).$ Since $\tau(\widetilde
S_0) = \{ \, x \in \overline {R_d(S)} \, : \, [x] \neq \phi \,
\},$ there exists $M \in \cM_d(\TT^n)$ such that $M \subseteq
\overline {R_d(S)}$ and $x \in M.$ Choose $O \subset \TT^n$ open
such that $x \in O$ and $E|_O \, : \, O \rightarrow \TT^n$ is
injective and construct $N \equiv E(O \cap M).$ Then $N \in
\cM_d(\TT^n),$ $N \subseteq \overline {R_d(S)},$ and $E(x) \in N,$
hence $E(x) \in \tau(\widetilde S_0)$ and this proves the first
assertion. The second and fourth assertions follows since
$\widetilde E (M_x) = N_x$ and the fourth assertion follows since
$\widetilde E \, : \, \widetilde S_0 \rightarrow \widetilde S_0$
is uniformly continuous. The assertions concerning $\cC(\widetilde
S_0)$ follow from Proposition \ref{pro:loja} and the surjectivity
of $E \, : \, \overline {R_d(S)} \rightarrow {R_d(S)}.$
\end{pf}
\begin{defin}
\label{defin:posexp}
Let $(X,\rho)$ be a compact metric space. A continuous (not
necessarily surjective) map $f \, : \, X \rightarrow X$ is called
expanding (with respect to $\rho$) if there exists numbers
$\epsilon > 0$ and $\lambda > 1$ such that $0 < \rho(x,y) <
\epsilon$ implies $\rho(f(x),f(y)) > \lambda \, \rho(x,y)$ and is
called positively expansive if there is a constant $c
> 0$ such that if $x, y \in X$ and $x \neq y$ then
$d(f^i(x),f^i(y))
> c$ for some $i \geq 0.$
\end{defin}
\begin{rem}
\label{rem:expanding}
Coven and Reddy \cite{covred} proved that a positively expansive
map of a closed (compact without boundary) topological manifold is
an expanding map with respect to some metric. Gromov \cite{gromov}
proved that an expanding differentiable map of a closed smooth
manifold is topologically conjugate to an expanding
infra-nil-endomorphism. Hiraide \cite{hiraide1} proved that a
positively expansive map of a closed topological manifold is
topologically conjugate to an expanding infra-nil-endomorphism.
\end{rem}
\begin{pro}
\label{pro:hiraide}
If $X$ is a metric space that is compact, connected, and locally
connected, and $f \, : \, X \rightarrow X$ is a positively
expansive map, and if $K \subset X$ is compact and satisfies the
following conditions:
\begin{enumerate}
\item $f(X \, \backslash \, K) \subseteq X \, \backslash \, K,$
\item $F|_{X \, \backslash \, K} \, : \, X \, \backslash \, K
\rightarrow X \, \backslash \, K$ is an open map,
\end{enumerate}
then $K = \phi.$
\end{pro}
\begin{pf}
This result was proved by Hiraide in \cite{hiraide2}, p. 566.
\end{pf}
\begin{thm}
\label{thm:ressing}
If $E \in \cE_n,$ $S \in \cS(\TT^n),$ and $E(S) = S$ then each $C
\in \cC(\widetilde S_0)$ is a compact manifold (without boundary),
$\tau(C) \in \cS(\TT^n)$ is irreducible, $S_0 = S,$ and
$\overline {R_d(S)} = \bigcup_{C \in \cC(\widetilde S)} \tau(C)$
is the union of immersed real analytic manifolds.
\end{thm}
\begin{pf}
We use Lemma \ref{lem:lift} to choose $p \in \NN$ such that
$\widetilde E^p(C) = C, \ C \in \cC(\widetilde S_0).$
Choose $C \in \cC(\widetilde S_0)$ and construct $X \equiv
\overline C$ and $K \equiv X \, \backslash \, C.$ Clearly $X$ is
compact, connected and locally connected metric space, the map $f
\equiv \widetilde E^p \, : \, X \rightarrow X$ is positively
expansive, $C = X \, \backslash \, K,$ and Brouwer's theorem on
invariance of domain implies that $f|_C \, : \, C \rightarrow C$
is open. Therefore $X,$ $f,$ and $K$ satisfy the hypothesis of
Proposition \ref{pro:hiraide} hence $K = \phi$ and $X = C.$ The
remaining assertions are obvious.
\end{pf}
\begin{rem}
\label{rem:resolution}
Grauert's theorem \cite{grauert}, Theorem 3, ensures that
$\widetilde S$ admits a closed analytic embedding in $\RR^m$ for
sufficiently large $m \in \NN.$ This fact can be used together
with a tubular neighborhood construction to provide a resolution
of singularities, in the sense of Hironaka \cite{krantz},
\cite{hironaka}, for the set $\overline {R_d(S)}.$ The existence
of the positively expansive map makes this resolution particularly
easy. Also see the discussion about embeddings into affine space
in \cite{andradas}, p. 224.
\end{rem}
We recall that the differential $d\tau$ maps tangent vectors of
$C$ into $\RR^n.$
\begin{cor}
\label{cor:ressing1}
Assume that $E \in \cE_n,$ $S \in \cS(\TT^n),$ $E(S) = S,$ $C \in
\cC(\widetilde S),$ ${\widetilde E}(C) = C,$ and for every $c \in
C$ at least one of the following conditions hold:
\begin{enumerate}
\item $d\tau(T_c) \in V_{\sigma}^{\perp},$
\item there exists $H_{c} \in \cG_c(\TT^n)$ such that
$\hbox{dim}(H_c) \geq 1$ and $\tau(x) + H_{c} \subseteq \tau(C).$
\end{enumerate}
Then, by induction on $n$ and $\hbox{dim}(S),$ $\tau(C) \in
\cF(\TT^n).$
\end{cor}
\begin{pf}
The hypothesis above implies that $C$ can be expressed as the
countable union
$$C = C_1 \cup \bigcup_{H \in \cG_1(\TT^n)}
\tau^{-1}(\tau(C)_H)$$
where $C_1$ denotes the subset of points $x$ in $C$ that satisfy
condition 1 above and where $\cG_1(\TT^n)$ denotes the set of
closed connected subgroups of $\TT^n$ whose dimension is $\geq 1.$
Since $C$ is a nonempty complete metric space, the Baire Category
Theorem implies that either $C_1$ or one of the sets
$\tau^{-1}(\tau(C)_H)$ has a nonempty interior. Since each of
these sets is a real analytic subset of the irreducible real
analytic manifold $C,$ either $C_1 = C$ or $\tau^{-1}(\tau(C)_H) =
C$ for some $H \in \cG_1(\TT^n).$ If $C_1 = C$ then a simple
argument shows that there exists $z \in \TT^n$ such that $\tau(C)
\subseteq \pi_{n}(V_{\sigma}^{\perp}) + z.$ Define $K \equiv
\pi_n((V_{\sigma}^{\perp})_{rat}).$ Since $\tau(C)$ is an immersed compact manifold
it satisfies the hypothesis of Proposition \ref{pro:irr}.
Therefore there exists $y \in \TT^n$ such that $E(y) = y$ and
$\tau(C)-y \subseteq K.$ Clearly $K \in \cG_c(\TT^n),$ and $K$ is
isomorphic to $\TT^m$ where $m < n,$ $E(K) = K,$ $\tau(C) - y \in
\cS(H),$ and $E(\tau(C)-y) = \tau(C)-y.$ By induction on $n$ it
follows that
$\tau(C)-y \in \cF(\TT^n)$
hence $\tau(C) \in \cF(\TT^n).$ If for some
$H \in \cG_1(\TT^n),$
$\tau^{-1}(\tau(C)_H) = C$
then $\tau(C)_H = \tau(C).$ Then Lemma \ref{lem:inv3} and
induction on $\hbox{dim}(S)$ imply that $\tau(C) \in \cF(\TT^n).$
\end{pf}
\begin{cor}
\label{cor:ressing2}
Assume that $E \in \cE_n,$ $S \in \cS(\TT^n),$ $E(S) = S,$ $C \in
\cC(\widetilde S),$ ${\widetilde E}(C) = C,$ $c \in C,$ $w \in
T_c(C),$ and $d\tau(w) \notin V_{\sigma}^{\perp}.$ Then there
exists $K \in \cG_c(\TT^n)$ such that $\hbox{dim}(K) \geq 1$ and
$\tau(c) + K \subset \tau(C).$
\end{cor}
\begin{pf}
Let $d = \hbox{dim}(S).$ The construction of $\widetilde S$
ensures that there exist $M \in \cM_d(\RR^n),$ $x \in M,$ and $v
\in T_x(M)$ such that $\pi_n(M) \subseteq \tau(C),$ $\pi_n(x) =
\tau(c),$ and $v = d\tau(w).$ From the argument used in
Proposition \ref{prop:asyminv} and the fact that $C$ is compact,
there exists a function $j \, : \, \NN \rightarrow \NN$ that
satisfies the following properties
\begin{enumerate}
\item $\lim_{i \rightarrow \infty} j(i) = \infty,$
\item there exists $y \in \tau(C)$ such that the sequence
where $y_i \equiv E^{j(i)}(\pi_n(x)),$ converges to $y,$
\item there exists $u \in \RR^n$ such that the sequence
$u_i \equiv \frac{E^{j(i)}v}{||E^{j(i)}v||}$ converges to $u.$
\item there exists $z \in C$ such that the sequence $c_i \equiv {\widetilde
E}^{j(i)}(c)$ converges to $z.$
\end{enumerate}
and there exists $H \in \cG_c(\TT^n)$ such that $\hbox{dim}(H)
\geq 1$ and $y + H \subseteq \tau(C).$ Clearly, since $\tau \circ
\widetilde E = E \circ \tau \, : \, C \rightarrow \tau(C),$
$\tau(c_i) = y_i$ converges to $y.$ Therefore, since $c_i$
converges to $z,$ $\tau(z) = y.$ Let $\rho$ be any metric on $C.$
Since $\tau$ is an \'{e}tale map, there exists a sequence $O_k
\subset C, \, k \in \NN$ of open neighborhoods of $c$ such that
the restrictions $\tau|_{O_k} \, : \, O_k \rightarrow \tau(O_k)$
are homeomorphisms and the $\rho$ diameters of $O_k$ converge to
zero. Since $\widetilde E \, : \, C \rightarrow C$ is expansive
and $C$ is a manifold, Brouwer's theorem on invariance of domain
implies that for each $k \in \NN$ there exists an integer $p(k)
\in \NN$ such that $z \in {\widetilde E}^{p(k)}(O_k).$ Therefore
$y \in E^{p(k)}(\tau(O_k)).$ For $k \in \NN,$ $(E^{p(k)})^{-1}(H)
\in \cG(\TT^n)$ and we let $H_k$ denote its connected component
that contains the identity. Therefore $H_k \in \cG_c(\TT^n)$ and
$\hbox{dim}(H_k) \geq 1.$ Since $\tau(C) \in \cS(\TT^n),$ there
exists $t_k \in \tau(O_k)$ such that $t_k + H_k \subseteq
\tau(C).$ Since the Hausdorff topology on $\cG_c(\TT^n)$ is
compact, there exists $K \in \cG_c(\TT^n)$ with $\hbox{dim}(K)
\geq 1$ and there exists a subsequence of $H_k$ that converges to
$K.$ Since $t_k$ converges to $\tau(c),$ $\tau(c) + K \subseteq
\tau(C).$
\end{pf}
We now prove Theorem \ref{thm:redmain}. It suffices to prove that
if $E$ and $S$ satisfy the hypothesis of the theorem then
$\overline{R_d(S)} \in \cF(\TT^n).$ Corollary \ref{cor:ressing1}
together with Corollary \ref{cor:ressing2} implies that $\tau(C)
\in \cF(\TT^n)$ whenever $C \in \cC(\widetilde S).$ Theorem
\ref{thm:ressing} implies that $\overline{R_d(S)}$ equals a finite
union of $\tau(C).$ The proof is complete.
\section{Appendix A. The Smith Normal Form.}
\setcounter{equation}{0}
For $m,n \in \NN$ we let $M_{n,m}(\ZZ)$ denote the set of $n$ by
$m$ integer matrices. A minor of order $k \in \NN$ of a matrix is
the determinant of a $k$ by $k$ submatrix and the elementary
divisor of order $k$ of a matrix, denoted by $d_k,$ is the
greatest common divisor of its minors of order $k.$ The
Cauchy-Binet theorem \cite{turnbull} implies that if $M \in
M_{m,n}(\ZZ)$ with $m \leq n$ then its elementary divisors satisfy
$d_1 | d_2 | \cdots | d_m.$ The matrix $M$ is called unimodular if
$d_m = 1.$ We let $U_{n,m}(\ZZ)$ denote the set of all unimodular
matrices in $M_{n,m}(\ZZ).$ Clearly $U_{n,n}(\ZZ)$ is the set of
all matrices in $M_{m,n}(\ZZ)$ whose determinant equals $\pm 1$
and it forms a group under matrix multiplication.
\begin{thm}
\label{thm:smith}
    If $m, n \in \NN, \, m \leq n,$ and $M \in M_{n,m},$
    then there exist $U_n \in U_{n,n}(\ZZ)$ and $U_m
    \in U_{m,m}(\ZZ)$ such that $U_n \, M \, U_m = D$
    where $D$ has the form
$$
D = \left[
\begin{array}{cccccc}
 d_1 & 0 &  & \cdots & 0 & 0
\\[1ex]
 0 & d_2 &  & \cdots & 0 & 0
\\[1ex]
\vdots & & & & & \vdots
\\[1ex]
0 & 0  &  & \cdots  & d_{m-1} & 0
\\[1ex]
0 & 0  &  & \cdots & 0 & d_m
\\[1ex]
0 & 0 &  & \cdots & 0 & 0
\\[1ex]
\vdots & & & & & \vdots
\\[1ex]
0 & 0 &  & \cdots & 0 & 0
\end{array}
\right]
$$
and $d_j, j = 1,...,m,$ are the elementary divisors of $M.$
\end{thm}
\begin{pf}
This decomposition, derived in 1861 by Smith \cite{smith}, is
described in \cite{newman}.
\end{pf}
\begin{cor}
\label{cor:GtoH}
If $G \in \cG_c(\TT^n)$ then there exists $H \in \cG_c(\TT^n)$
such that $G \cap H = \{0\}$ and $G + H = \TT^n.$
\end{cor}
\begin{pf}
Construct $V \equiv \pi_{n}^{-1}(G).$ Then since $\pi_n(V) = G \in
\cG(\TT^n)$ Lemma \ref{lem:rational} implies that $V$ is a
rational subspace of $\RR^n,$ hence $V$ is spanned by vectors in
$V \cap \ZZ^n.$ Let $m \equiv dim(V)$ and let $M \in M_{n,m}$ such
that the columns of $M$ span $V.$ Theorem \ref{thm:smith} implies
that there exists $U_n \in U_{n,n}(\ZZ)$ and $U_m \in
U_{m,m}(\ZZ)$ such that $U_n \, M \, U_m = D$ where $D \in
M_{n,m}(\ZZ)$ is the diagonal matrix described in Theorem
\ref{thm:smith}. Therefore, since $V$ is spanned by the columns of
$U_{n}^{-1}DU_{m}^{-1},$ $V$ is spanned by the first $m$ columns
of the matrix $U_{n}^{-1}.$ Let $W$ denote the subspace $W$
spanned by last $n-m$ columns of $U_{n}^{-1}$ and let $H =
\pi_n(W).$ Since $W$ is a rational subspace, Lemma
\ref{lem:rational} implies that $H \in \cG_c(\TT^n).$ Since $V
\cap W = \{0\}$ and $V + W = \TT^n$ it follows that $H$ satisfies
the properties asserted above.
\end{pf}
\section{Appendix B: Fourier Analysis and Hausdorff Topology on
$\cG_c(\TT^n).$}
\setcounter{equation}{0}
In this appendix we prove Theorem \ref{thm:invinv} by exploiting a
compact Hausdorff topological space structure on $\cG_c(\TT^n)$
that we construct using Pontryagin duality theory. This theory was
initially developed by Lev Semenovich Pontryagin
\cite{pontryagin1}, \cite{pontryagin2} to extend the classical
Fourier analysis to compact and discrete abelian groups. It was
extended later to locally compact abelian groups by E. R. van
Kampen \cite{kampen1}, \cite{kampen2} and others. If $E \in
M_{n,n}(\ZZ)$ then $E$ induces maps on $\cG_c(\TT^n), on
\cG(\TT^n),$ and on $\cF(\TT^n).$ If $E \in U_{n,n}(\ZZ)$ then
$E^{-1}$ induces maps on $\cG(\TT^n)$ and on $\cF(\TT^n).$ We
record the following standard results \cite{hewitt}, \cite{rudin}.
\begin{enumerate}
\item We let $\UU$ denote the multiplicative group of complex
numbers of modulus one. For $G \in \cG(\TT^n)$ its Pontryagin dual
$G^{\wedge}$ consists of all continuous homomorphisms $\chi \, : \,
G \rightarrow \UU$ under pointwise multiplication and equipped with
the discrete topology.
\item The map $\ZZ^n \ni \ell \rightarrow \chi_{\ell} \in (\TT^n)^{\wedge}$
defined by
$\chi_{\ell}(x + \ZZ^n) \equiv \hbox{exp}(2\, \pi \, i \, \ell \cdot
x), \ x \in \RR^n$
is an isomorphism. We will identify $(\TT^n)^{\wedge}$ with $\ZZ^n.$
\item If $G \in \cG(\TT^n)$ then the dual $i_{G}^{\wedge} \, : \ZZ^n \rightarrow
G^{\wedge},$ defined by
$
    i_{G}^{\wedge}(\chi) \equiv \chi|_{G}, \ \ \chi \in \ZZ^n,
$
of the inclusion map $i_{G} \, : \, G \rightarrow \TT^n,$ is an
epimorphism. We let $G^{\perp}$ denote the kernel of
$i_{G}^{\wedge}.$ The map $G \rightarrow G^{\perp}$ is a bijection
between $\cG(\TT^n)$ and the set of subgroups of $\ZZ^n.$ Therefore
$\cG(\TT^n)$ is countable.
\item If $G \in \cG(\TT^n)$ then $G^{\wedge}$ is
isomorphic to $\ZZ^n \, / \, G^{\perp}.$ This is a finitely
generated abelian group that is isomorphic to $\ZZ^{n-m} \oplus H,$
where $m$ is the rank of $G^{\perp}$ and $H$ is a finite group.
\item Let $G$ and $H$ be the groups above and let $G_c$ denote the
connected component of $G$ that contains the identity. The group
$G_{c}^{\perp}$ is isomorphic to $\ZZ^m,$ the group $G_c$ is
isomorphic to $\TT^{m},$ and the group $H$ is isomorphic to each
of the following groups: $G \, / \, G_c,\ $ $G_{c}^{\perp} \, / \,
G^{\perp},$ $\pi_{n}^{-1}(G) \, / \, \pi_{n}^{-1}(G_c).$
\end{enumerate}
\begin{lem}
\label{lem:FinV}
$\cF(\TT^n) \subset \cV(\TT^n).$
\end{lem}
\begin{pf}
We first show that if $G \in \cG_c(\TT^n)$ and $y \in \TT^n,$ then
$G + y \in \cV_{irr}(\TT^n).$ Let $k = \hbox{dim}(G).$ Then
$G^{\perp}$ is a rank $n-k$ subgroup of $\ZZ^n.$ Choose
$\chi_1,...,\chi_{n-k} \in \ZZ^n$ that generate $G^{\perp}.$ We
observe that $G = \{ x \in \TT^n \, : \, \chi_j(x) = 1, \ j =
1,....,n-k \, \}.$ The function $h \, : \, \TT^n \rightarrow \CC$
defined by
$ h(x) \equiv \sum_{j=1}^{\, n-k} \, | \, \chi_{j}(x) -
\chi_{j}(y) \, |^2 $
is in $\cA(\TT^n)$ and $G + y = \cZ_h.$ Therefore $G+y \in
\cV(\TT^n).$ If $a,b \in \cA(\TT^n)$ and $G+y \subseteq \cZ_{a}
\cup \cZ_{b}$ then either $a$ or $b$ must vanish on a subset of
$G+y$ that is open (relative to the topology on $G+y$ induced as a
subset of $\TT^n$). Since both $a$ and $b$ are real analytic one
of them must vanishes on $G+y$ and this shows that $G+y \in
\cV_{irr}(\TT^n).$ Therefore $\cF(\TT^n) \subset \cV(\TT^n)$ since
$\cV(\TT^n)$ is closed under finite unions.
\end{pf}
We state the following result, which follows directly from Theorem
\ref{thm:smith}, without proof.
\begin{cor}
\label{cor:perp}
Let $G \in \cG(\TT^n)$ and let $G_c \in \cG_c(\TT^n)$ denote the
connected component of $G$ that contains the identity. Define the
subspace $V \equiv \pi_{n}^{-1}(G_c)$ of $\RR^n$ and define the
lattice subgroup $L \equiv V \cap \ZZ^n$ of $\ZZ^n.$ Let $m$
denote the rank of $G^{\perp}$ and let $M \in M_{n,m}(\ZZ)$ such
that the columns of $M$ generate $G^{\perp}.$ Choose $U_n \in
U_{n,n}(\ZZ)$ and $U_m \in U_{m,m}(\ZZ)$ such that $U_n \, M \,
U_m = D$ where $D$ has the form in Theorem \ref{thm:smith}. Then
$G^{\perp}$ is spanned by the columns of the matrix $U_{n}^{-1} \,
D,$ $G_{c}^{\perp}$ is generated by the first $m$ columns of the
matrix $U_{n}^{-1},$ and $L$ is generated (and $V$ is spanned) by
the last $n-m$ columns of the matrix $U_{n}^{T}$ (the transpose of
$U_{n}$).
\end{cor}
We let $\cU$ denote the set of all subsets of $\cG(\TT^n) \times
\cG(\TT^n)$ that contain as a subset a set having the form
$
    O(F) \equiv \{ \, (G_1,G_2) \in \cG(\TT^n) \times \cG(\TT^n)
    \, : \, G_{1}^{\perp} \cap F = G_{2}^{\perp} \cap F \, \}
$
where $F \subset \ZZ^n$ is finite. We observe that
$(\cG(\TT^n),\cU)$ is a uniform space that is Hausdorff and
metrizable (\cite{kelly}, Chapter 6).
\begin{lem}
\label{lem:Cauchy}
A sequence $G_j \in \cG(\TT^n), \ j \in \NN$ is a Cauchy sequence
if and only if for every finite $F \subset \ZZ^n,$ there exists
$J(F) \in \NN$ such that $(G_j,G_k) \in O(F), \ \ j,k \geq J(F).$
If $G_j \in \cG(\TT^n)$ is a Cauchy sequence then the following
properties hold:
\begin{enumerate}
\item the set
$ C \equiv \bigcup_{k \geq 1} \, \bigcap_{j=k}^{\infty}
G_{j}^{\perp}$
is a subgroup of $\ZZ^n,$
\item $G_j$ converges to $G \in \cG(\TT^n)$ defined by: $G^{\perp} =
C,$
\item if $F \subset \ZZ^n$ is a finite set that generates
$G^{\perp}$ and $j \geq J(F)$ then $G_j \subseteq G,$
\item if $G_j \in \cG_c(\TT^n), \ j \in \NN,$ then $G \in \cG_c(\TT^n).$
\end{enumerate}
\end{lem}
\begin{pf}
The first assertion is obvious. The second assertion follows since
if $F \subset \ZZ^n$ is finite and if $j \geq J(F)$ then $(G_j,G)
\in O(F).$ The third assertion follows since if $F \subset \ZZ^n$
is finite and generates $G^{\perp}$ and $j \geq J(F)$ then $F =
G^{\perp} \cap F = G_{j}^{\perp} \cap F$ hence $F \subset
G_{j}^{\perp}$ hence $G^{\perp} \subseteq G_{j}^{\perp}$ hence
$G_{j} \subseteq G.$ The fourth assertion follows since if $G_j$
is connected and $G_j \subseteq G$ then $G_j \subseteq G_c.$
\end{pf}
\begin{thm}
\label{thm:compact}
The uniform spaces $(\cG(\TT^n),\cU)$ and $(\cG_c(\TT^n),\cU)$ are
compact.
\end{thm}
\begin{pf}
Let $\cH(\TT^n)$ denote the set of all closed subsets of $\TT^n,$
let $\rho$ be any metric on $\TT^d$ that induces the standard
topology on $\TT^n,$ and define the associated Hausdorff metric on
$\cH(\TT^n)$ as follows
\begin{equation}
\label{hausdorff}
    \rho_H(A,B) \equiv
    \hbox{max} \{ \, \sup_{x \in A} \, \inf_{y \in B} \rho(x,y),
    \, \sup_{y \in B} \, \inf_{x \in A}
    \rho(x,y) \, \}, \ \ A, B \in \cH(\TT^n).
\end{equation}
The metric space $(\cH(\TT^n),\rho_H)$ is compact
(\cite{dugundji}, p.205, p.253), (\cite{munkres}, p.279).
Furthermore, the uniformity on $\cH(\TT^n)$ defined by $\rho_H$ is
independent of the metric $\rho$ on $\TT^n$ and its restriction to
$\cG(\TT^n)$ coincides $\cU.$ Since $(\cG(\TT^n),\cU),$ and
$(\cG_c(\TT^n),\cU)$ are complete, they are closed and therefore
compact.
\end{pf}
\begin{rem}
Since $\TT^n$ is a compact Hausdorff space, the Hausdorff topology
on the space $\cH(\TT^n)$ of all closed subsets of $\TT^n$
coincides with both the Vietoris topology and the Fell topology
\cite{fell}, \cite{michael}. We will denote the topological space
$(\cG(\TT^n),\cU), (\cG_c(\TT^n),\cU)$ simply by $\cG(\TT^n),
\cG_c(\TT^n)$ respectively.
\end{rem}
\begin{thm}
\label{thm:invinv}
If $E \in \cE_n,$ $S \subseteq \TT^n,$ $\overline S = S,$ $E(S) =
S,$ $H \in \cG_c(\TT^n),$ $dim(H) \geq 1,$ and $S = S_H,$ then there
exists $G \in \cG_c(\TT^n)$ and $p \in \NN$ such that $dim(G) \geq
1,$ $E^p(G) = G,$ and $S = S_G.$
\end{thm}
\begin{pf}
We first observe that $E$ induces an injection $E \, : \,
\cG_c(\TT^n) \rightarrow \cG_c(\TT^n).$ Construct the orbit $X
\equiv \{ \, E^j(H) \, : \, j \geq 0 \, \}$ of $H$ under $E$ and
let $\overline X$ denote the topological closure of $X$ in
$\cG_c(\TT^n).$ Clearly $E(X) \subseteq X.$  Choose $G \in
\overline X$ so that $dim(G) \geq dim(D)$ for all $D \in \overline
X$ and define $Y \equiv \{ \, E^j(G) \, : \, j \geq 0 \, \}.$
Property (3) in Lemma \ref{lem:Cauchy} implies that if $Y$ is
infinite then there would exists $D \in \overline Y$ with $dim(D)
> dim(G)$ and this contradicts the choice of $G.$ Therefore $Y$ is
finite. Since $E(Y) \subseteq Y$ the map $E|_{Y} \, : \, Y
\rightarrow Y$ is an injection and therefore a bijection.
Therefore there exists $p \in \NN$ such that $E^p(G) = G.$ Let
$n_j \in \NN$ be a sequence such that $E^{n_j}(H) \rightarrow G.$
Then $S = S_H$ implies $S = S + H$ hence $S = E^{n_j}(S) = S +
E^{n_j} \rightarrow S + G$ hence $S = S_G.$
\end{pf}
\section{Appendix C: Lojasiewicz's Structure Theorem}
\setcounter{equation}{0}
We state, with slight modifications, and derive consequences of
Lojasiewicz's structure theorem for real analytic varieties
\cite{lojaciewicz}, as presented by Krantz in \cite{krantz}, p.
152-156. If $k \geq 1$ and $U \subseteq \RR^{k-1}$ is open (note
that $\RR^{0} \equiv \{0\}$ hence $\cA(\RR^0) = \RR$), a function
$H \, : \, U \times \RR \rightarrow \RR$ is called a {\it
distinguished polynomial defined on} $U$ if there exist $m \in
\NN,$ $c_j \in \cA(U), \ 0 \leq j \leq m-1$ such that
$$
    H(x,y) = y^m + \sum_{j=0}^{m-1} c_{j}(x)
    y^{j}, \ \ x \in U, \,
    y \in \RR.
$$
\begin{pro}
\label{pro:weierstrass}
(Weierstrass Preparation Theorem) If $f$ is real analytic in a
open neighborhood of $0 \in \RR^k$ and $f(0,\cdots,0,x_k) \neq 0$
then there exists an open neighborhood $U$ of $0 \in \RR^{k-1},$ a
positive number $\delta,$ a function (called a unit) $u \in \cA(U
\times (-\delta,\delta))$ that never vanishes, and a distinguished
polynomial $H$ defined on $U$ such that
$$
    f(x,y) = u(x,y) \, H(x,y), \ \ x \in U, \, y \in
    (-\delta,\delta).
$$
\end{pro}
We observe that since the ring of (germs of) real analytic
functions in a neighborhood of $0$ is a unique factorization
domain, the ring of distinguished polynomials in a neighborhood of
$0$ is also a unique factorization domain. Furthermore, a
distinguished polynomial $H$ defined on $U$ has repeated factors
iff its discriminant vanishes (everywhere) on $U.$ Therefore,
given a distinguished polynomial $H$ defined on $U,$ there exists
a unique polynomial $H_0$ defined on $U$ whose discriminant does
not vanish everywhere and whose zeros coincide with the zeros of
$H.$ The following renown result provides a detailed description
of the relationship between real analytic subsets and real
analytic submanifolds:
\begin{pro}
\label{pro:loja} (Lojaciewicz's Structure Theorem for Real
Analytic Varieties) Let $O$ be an open neighborhood of $0 \in
\RR^n$ and let $f \in \cA(O)$ satisfy $f(0,\cdots,0,x_{n-1}) \neq
0$ and $d \equiv dim(\cZ_f,0) < n.$ If $d = 0$ then $\cZ_f$ is a
discrete set. If $d \geq 1$ then after a suitable rotation of
$\RR^n$ that effects only the the first $n-1$ coordinates, there
exist positive numbers
$
    \delta_1,\cdots,\delta_n
$
and for every $1 \leq k \leq d,$ $U_k = \prod_{\, i=1}^{\, k}
(-\delta_i,\delta_i),$ a system of distinguished polynomials
$
    H_{\ell}^{k}(x,y), \ \ x \in U_k, \, y \in \RR, \, k+1 \leq \ell \leq n $
defined on $U_k$ that satisfy the following properties:
\begin{enumerate}
\item $f \in \cA(U_n),$
\item the discriminant of $H_{\ell}^{k}$ does not vanish everywhere on $U_k,$
\item if $H_{\ell}^{k}(x,y) = 0$ and $x \in U_k$ then $y \in
(-\delta_{\ell},\delta_{\ell}),$
\item there exists pairwise disjoint $V^k \in \cM(\RR^n), \ k = 0,\cdots,d$ such that
$
    U_{n} \cap \cZ_f = V^0 \cup V^1 \cup \cdots \cup V^d,
$
\item either $V^0 = \phi$ or $V^0 = \{0\},$
\item for $1 \leq k \leq d$ there exists $N_k \in \NN$ and pairwise disjoint
connected $k$ dimensional real analytic submanifolds
$
    \Gamma_{j}^{k} \subset \RR^n, \ \ j = 1,\cdots,N_k
$
such that
$
    V^k = \bigcup_{\, j = 1}^{\, N_k} \Gamma_{j}^{k},
$
\item (Analytic Parameterization) For $1 \leq k \leq d$ and $1 \leq j \leq N_k$ there
exist a open set $U_{j,k} \subseteq U_k$ and real analytic
functions $\gamma_{\ell,j,k} \in \cA(U_{j,k}), \  \ k+1 \leq \ell
\leq n$ such that
$
    \Gamma_{j}^{k} = \{
    (x,\gamma_{k+1,j,k}(x),\cdots,\gamma_{n,j,k}(x)) \, :
    \, x \in U_{j,k} \, \},
$
\item $H_{\ell}^{k}(x,\gamma_{\ell,j,k}(x)) = 0, \ \ x \in U_{j,k},$
\item the discriminant of $H_{\ell}^{k}(x,y)$ satisfies $D_{\ell}^{k}(x) \neq 0, \  \ x \in U_{j,k},$
\item (Non-Redundancy) For $1 \leq k \leq d$ and $1 \leq i < j \leq
N_k$ either $U_{i,k} \cap U_{j,k} = \phi$ or $U_{i,k} = U_{j,k}.$
In the latter case for $k+1 \leq \ell \leq n$ either
$\gamma_{\ell,i,k} = \gamma_{\ell,j,k}$ or $\gamma_{\ell,i,k}(x)
\neq \gamma_{\ell,j,k}(x), \ x \in U_{\ell,i,k}.$
\item (Stratification) For $1 \leq k \leq d$ and $1 \leq j \leq N_k$
the $U \cap \partial \Gamma_{j}^{k}$ is a union of sets of the
form $\Gamma_{i}^{p}$ with $1 \leq p < k$ and $1 \leq i \leq N_p$
and possibly $V^0.$
\end{enumerate}
\end{pro}
\begin{cor}
\label{cor:loja1}
If $S \in \cS(\TT^n),$ $d = dim(S),$ $x \in R_d(S),$ $M, N \in
\cM(\TT^n),$ $M \subseteq S,$ $N \subseteq S,$ and $x \in M \cap N,$
then either $M_x = N_x$ or there exists open $O \subset \RR^n$ such
that for every $y \in O \cap M \cap N,$ $(O \cap M)_y \neq (O \cap
N)_y.$
\end{cor}
\begin{pf}
Lojasiewicz's Structure Theorem implies that there exists open
neighborhood $O$ of $x$ such that if there exists $y \in O \cap M
\cap N$ such that $(O \cap M)_y = (O \cap N)_y$ then there exists
$\Gamma_{j}^{d}$ such that $O \cap \Gamma_{j}^{d} \subseteq O \cap M
\cap N$ and $x \in \overline \Gamma_{j}^{d}.$ The principle of
analytic continuation then implies that $M_x = N_x.$
\end{pf}
{\it Acknowledgments} The author thanks Zbigniew Jelonek for
informing him about the result of Frisch mentioned in Remark
\ref{rem:frisch}, and his colleague Zhang De-Qi for informing him
about the use of \'{e}tale maps in lifting constructions.

\end{document}